\documentclass[12pt,a4paper]{amsart}

\usepackage{amssymb}
\usepackage[final]{showkeys}
\usepackage{microtype}
\usepackage{color}
\usepackage{graphicx}
\usepackage[hmargin=3.5cm,vmargin={3.5cm,4cm}]{geometry}

\newtheorem{te}{Theorem}

\newtheorem{os}[te]{Remark}
\newtheorem{prop}[te]{Proposition}

\allowdisplaybreaks

\begin{document}

\title[FBm rules by nonlinear equations]{Fractional Brownian motions 
ruled by nonlinear equations}

\author[R. Garra]{Roberto Garra$^{1,*}$}
\address{{}$^1$Department of Statistical Science, Sapienza University of Rome, Italy}
\email{roberto.garra@sbai.uniroma1.it}

\author[E. Issoglio]{Elena Issoglio$^2$}
\address{{}$^2$School of Mathematics, University of Leeds, Leeds, UK} \email{e.issoglio@leeds.ac.uk}

\author[G.S. Taverna]{Giorgio S. Taverna$^3$}
\address{{}$^3$School of Earth and Environment, University of Leeds, Leeds, UK}
\email{eegst@leeds.ac.uk}

 \address{{}$^*$corresponding author}
 
    \date{\today}

    \begin{abstract}
	In this note we consider generalized diffusion equations in which the diffusivity coefficient is not necessarily constant in time, but instead it solves a nonlinear fractional differential equation involving fractional Riemann-Liouville time-derivative.
Our main contribution is to highlight the link between these generalised equations and fractional Brownian motion (fBm). In particular, we investigate  the governing equation of fBm and show that its diffusion coefficient must satisfy an additive evolutive  fractional equation. We derive in a similar way the governing equation of the iterated fractional Brownian motion.
\end{abstract}

	    
\maketitle

\section{Introduction}

Fractional Brownian motion (fBm in short) is an extension of classical Brownian motion (Bm), whose   formulation dates back to seminal work of Mandelbrot and Van Ness in 1968 \cite{Mandelbrot1968}. Subsequently, fBm and relative stochastic calculus were studied extensively, motivated by applications in numerous research areas and technological fields, see e.g.~\cite{biagini_et.al2008, ustunel, taqqu}.

In the last 25 years Fractional Calculus (FC) has been object of intensive research efforts with results spanning from economics, mathematics, physics, biology and chemistry, see e.g.~\cite{AnatolyA.KilbasHariM.Srivastava2006}.  
Although fBm and FC are independent areas of mathematics, they are strongly linked. One fundamental example of this interplay is the well known representation for fBm as fractional integral, which is one of the main tools to define stochastic integration with respect to fBm, see e.g.~the recent book \cite{biagini_et.al2008}. 

In this note we are particularly interested in recent links that have been discovered in mathematical physics and that  involve equations with time-dependent diffusivity coefficients \cite{Garra2015}. Following this idea, here we draw a new link between fractional Brownian motions and fractional differential equations. The starting point in our analysis is a \emph{generalized heat equation with varying diffusivity}. The idea of time-varying diffusivity is not new and has been successfully applied, for instance, in the study of mobility of biological systems \cite{A, B} and in the diffusion of pollutants in the atmospheric planetary boundary layer \cite{C}.
Our contribution here is to describe the  time-dependence of the diffusion coefficient by an \emph{additive nonlinear fractional equation}.  The coupling of these (generalised heat and nonlinear fractional) equations  leads to a system of governing equations for fractional Brownian motion. A similar approach is developed in order to obtain the governing equations of the iterated fractional Brownian motion.

The paper is structured as follows: 
In Section 2 we introduce  fBm and some of the basic properties and  equations related to it, while in Section 3 we illustrate the results obtained when the diffusivity term is the solution of a nonlinear fractional differential equation. In Section 4 we briefly conclude.

\section{Preliminaries}

It is well-known that the probability density function of  Brownian motion coincides with the fundamental solution to the heat equation. In what follows we recall the details of this fact for a slightly more general process, namely for a rescaled Bm $Y$ (which features a more general diffusion coefficient). The process $Y:=\{Y(t), t\geq 0\}$ is defined for all $t\geq 0$  by $Y(t):=\sqrt{2C}B(t)$, where $C$ is a positive constant and $\{B(t), t\geq 0\}$ is a standard Bm. 
The coefficient $\sqrt{2C}$  is called \emph{diffusion coefficient} because it is very much linked to the diffusivity of the related partial differential equation, as we will see below. 
When $C=\frac12$ we recover the standard Brownian motion case.
If we denote by $p(t,x)$ the density of the rescaled Brownian motion $Y(t)$, then $p$ satisfies the heat equation 
\begin{equation*}
\frac{\partial p}{\partial t} = C\frac{\partial^2 p}{\partial x^2},
\end{equation*}
and its explicit expression is 
\[
p(t,x) = \frac{1}{\sqrt{4  \pi C t} } e^{-\frac{|x|^2}{4Ct}}.
\]
For this process the mean-square displacement is linear in time, in particular one gets $E[Y(t)^2] = 2Ct$.

In this communication we are interested in a generalisation of the rescaled Brownian motion, namely  fractional Brownian motion. A fBm $B^H:=\{B^H(t), t\geq0\}$ is a family stochastic processes indexed by the  Hurst index $H\in(0,1)$. For each $H$ the process $B^H$ is defined to be a continuous and centered Gaussian process with covariance function $  E[B^H(t) B^H(s)] = \frac12 \left( t^{2H} + s^{2H} - |t-s|^{2H}\right)$ for all $t,s\geq0$. In the special case when $H=\frac12$ one recovers the classical Brownian motion $B$, with covariance 
 $  E[B (t) B (s)] = \frac12 \left( t  + s - |t-s| \right) = \min\{t, s\}.$
FBm is not a martingale nor a Markov process (unless $H=\frac12$) but it has the following interesting properties which do not hold for Bm: 
\begin{itemize}
\item (\emph{long-range dependence}) If $H>\frac12$ then the increments of fBm are positively correlated, if $H<\frac12$ increments are negatively correlated (and if $H=\frac12$ increments are independent). This produces effects of memory and persistence if $H>\frac12$ and of intermittency and anti-persistence if $H<\frac12$.

This memory feature allows fBm to be applied in many fields, for example telecommunication networks, weather derivatives in finance or filtering, see e.g.~\cite{ustunel}.
\item (\emph{power law mean-square displacement}) The variance (also known as mean-square displacement in Physics) at time $t$ is $E[(B^H(t))^2] = t^{2H}$. Equivalently the standard deviation is $t^{H}$. This means that, as the Hurst parameter changes, fBm can exhibit a power-law  mean-square displacement behaviour, hence nonlinear in $t$. In particular we can see powers of $t$ both bigger and smaller than 1, confirming the qualitatively different behaviour of fBm for $H>\frac12$ and $H<\frac12$. 

This feature is useful when studying phenomena whose mean-square displacement does not behave linearly in $t$. One example is the study of atmospheric dispersion of pollutant in an anomalous diffusive, due to the presence of turbulence, see  \cite{Goulart_et.al.2017}, where it is observed that data exhibits a non-linear mean-square displacement of the form $const\cdot t^\alpha$. 
\end{itemize}
The partial differential equation (PDE) related to fractional Brownian motion is known and it is a generalised heat equation with a time-dependent diffusion coefficient. As above we consider the rescaled process $Y^H(t) :=\sqrt{2C}B^H(t)$ to be slightly more general.
If we denote by $c_H(t, x)$ the density of the rescaled  fBm $Y^H$ with Hurst parameter $H$, then $c_H$  satisfies 
\begin{equation}\label{eq: heat eq c}
\frac{\partial c_H}{\partial t} = 2H C t^{2H-1}  \frac{\partial^2 c_H}{\partial x^2},
\end{equation}
and its explicit expression is 
\[
c_H(t,x) = \frac{1}{\sqrt{4\pi C  t^{2H}} } e^{-\frac{|x|^2}{4Ct^{2H}}}.
\] 
As expected, for this process the mean-square displacement is not linear in time (unless $H=\frac12$ which corresponds to a  rescaled  Brownian motion). In particular one gets $E[Y^H(t)^2] = 2Ct^{2H}$. 
When $C=\frac12$ this reduces to standard fBm. 

\section{The main result}

\subsection*{Fractional Brownian motion ruled by a nonlinear fractional differential equation}
Let us now generalise equation \eqref{eq: heat eq c} and consider the diffusion equation with a general time-varying diffusivity $D(t)$ in place of $ 2H C t^{2H-1} $, that is
\begin{equation}\label{1}
\frac{\partial \rho}{\partial t} = D(t)\frac{\partial^2 \rho}{\partial x^2}.
\end{equation}
Our aim is to reconstruct the dynamics \eqref{eq: heat eq c} for $c_H(t,x)$ that corresponds to fBm by means of this equation, in particular characterising it in terms of $\rho(t,x)$ given by \eqref{1} together with a dynamics for the coefficient $D(\cdot)$. It turns out (as we show below) that for $0<H<\frac12$ the correct assumption is that the time-dependence of dynamical diffusivity coefficient $D(\cdot)$ is ruled by the following nonlinear fractional differential equation
\begin{equation}\label{nfod}
\frac{d^{1-2H}D}{d t^{1-2H}} =k D^2, 
\end{equation}
 with 
\begin{equation}\label{4}
k = \frac{\Gamma(2H)}{ 2HC\Gamma(4H-1)}, \quad H\neq \frac14.
\end{equation}
Here $\Gamma(\cdot)$ denotes the Gamma function. 
Equation \eqref{nfod} involves a time-fractional derivative in the sense of Riemann-Liouville, i.e.
\begin{equation*}
\frac{d^{1-2H}D(t)}{d t^{1-2H}} = \frac{1}{\Gamma(2H)}\frac{d}{dt}\int_0^t\frac{D(s)ds}{(t-s)^{1-2H}},
\end{equation*}
for $ 0<H<\frac12$. 
Observing that
\begin{equation*}
\frac{d^{1-2H}t^\beta}{d t^{1-2H}} = \frac{\Gamma(\beta+1)}{\Gamma(\beta+2H)}t^{\beta-1+2H},
\end{equation*}
the solution of \eqref{nfod} is given by 
\begin{equation*}
D(t) = 2HC t^{2H-1}.
\end{equation*}
We note that the real order of the time-fractional nonlinear ordinary differential equation \eqref{nfod} is of course related to the Hurst index, in particular we take a fractional derivative of order $1-2H$ (which is positive for $0<H<\frac12$).

We summarise these observations in the proposition below. Note that for $C=\frac12$ we recover standard fBm. 
\begin{prop} The governing equation of the rescaled  fractional Brownian motion $Y^H$ for $H\in(0,\frac12)\setminus \{\frac14\}$ is given by \eqref{1}, where the diffusion coefficient $D(\cdot)$ satisfies the fractional differential equation \eqref{nfod}.
\end{prop}
By means of this kind of decomposition, we have a fractional Brownian motion ruled by a nonlinear fractional equation that provides the dynamical condition that must be satisfied by the diffusivity coefficient in order to obtain the diffusion equation governing a fractional Brownian motion. This connection between fractional Brownian motions and fractional derivative equations is new and can give new insights. The diffusivity coefficient is ruled by a nonlinear fractional differential equation in the case of negatively correlated increments $0<H<\frac12$ (excluding the value $H = \frac14$ that leads to the divergence of the Gamma function in \eqref{4}). Below we will see in Proposition \ref{pr 2} that a similar connection holds    in the case of positively correlated increments $\frac12< H<1$ by means of fractional integral equations. 

\begin{os} 
We now comment on the extreme case $H=0$ (and $C=\frac12$), which has been  excluded here so far. 
We first observe that
\begin{equation*}
k =  \frac{2H\Gamma(2H)}{4H^2C\Gamma(4H-1)} = \frac{\Gamma(2H+1)\cdot(4H-1)}{HC \cdot\Gamma(4H+1)}
\end{equation*}
and therefore, for $H\rightarrow 0$, we have that $k \rightarrow -\infty$ and equation \eqref{nfod} would lose its meaning. 
This is   in line with the fact that the classical definition of fBm via covariance is meaningless in the extreme case of  $H=0$, because it leads to a Gaussian process with constant  covariation $  E[B^H(t) B^H(s)] =1 $.
\end{os}

Next we want to obtain the governing equation of  fractional Brownian motion in the case of positively correlated increments, namely for  $\frac12<H<1$. To this aim we observe that now the quantity  $1-2H$ is negative and so it is natural to expect that the differential  equation governing the diffusion coefficient $D$ will involve a  Riemann-Liouville \emph{integral} of positive order $2H-1 $. This integral is defined by
\begin{equation*}
J_t^{2H-1}D(t) = \frac{1}{\Gamma(2H-1)}\int_0^t \frac{D(s)ds}{(t-s)^{2H-2}}.
\end{equation*}
Analogously as the fractional derivative case, we have that the diffusivity coefficient $D(\cdot)$ must satisfy the following nonlinear integral equation
\begin{equation}\label{2}
J_t^{2H-1}D(t) =k D^2 (t), 
\end{equation}
with $k $ as in \eqref{4}.
Recalling that, for $\beta>-1$
\begin{equation*}
J_t^{2H-1} t^\beta = \frac{\Gamma(\beta+1)}{\Gamma(\beta+2H-1)}t^{\beta+2H-1},
\end{equation*}
one obtains that the solution is given by  $D(t) = 2HC t^{2H-1}$. 
Summarising the result we have the proposition below. 
\begin{prop} \label{pr 2}
The governing equation of the rescaled  fractional Brownian motion $Y^H$ for $H\in(\frac12,1)$ is given by \eqref{1}, where the diffusion coefficient $D(\cdot)$ satisfies the fractional integral equation \eqref{2}.
\end{prop}

\begin{os}
Again we can look at the extreme case $H=1$ (and $C=\frac12$). In this case it is known that fBm $B^1$ with the given covariance $E[B^1(t) B^1(s)] = \frac12 \left( t^{2} + s^{2} - |t-s|^{2}\right) $ can be constructed by $B^1(t) = t\xi$ where $\xi$ is a standard Gaussian random variable. If we look at our representation in terms of  $D(\cdot)$ we see that $D$ must satisfy a classical integral equation (no fractional integral) of the form $\int_0^t D(s) ds =k D^2 (t)$. The Gamma functions appearing in the constant $k$ are finite and the solution $D(t)$ is linear  in $t$.  
\end{os}

\subsection*{Iterated Fractional Brownian motion ruled by nonlinear partial differential equation}
In the recent mathematical literature many papers have been devoted to the interplay between iterated processes and the governing partial differential equations, starting from the introduction of the iterated Brownian motions (we refer for example to \cite{DeBlassie} and references therein).
In   \cite{beghin} and \cite{mirko}, interesting connections have been pointed out between various types of compositions involving fractional Brownian motions and partial differential equations. 
In particular, it was shown that the distribution
\begin{equation*}
\rho(t,x) = 2\int_{0}^{\infty}\displaystyle\frac{e^{-\frac{x^2}{2s^{2H_1}}}}{\sqrt{2\pi s^{2H_1}}}\frac{e^{-\frac{s^2}{2t^{2H_1}}}}{\sqrt{2\pi t^{2H_2}}}ds
\end{equation*}
 of the iterated fractional Brownian motion $B_1^{H_1}(|B^{H_2}_2(t)|)$ is a solution of the first order
partial differential equation
\begin{equation*}
\frac{\partial \rho}{\partial t} = -H_1 H_2 \frac{\partial}{\partial x}\left(\frac{x}{t}\rho\right), \quad x\in \mathbb{R}, \ t>0, \ H_1, H_2 \in (0,1).
\end{equation*}
In analogy with the previous discussion, we can recast this equation
as a system of coupled equations, with a space-time varying  diffusion coefficient governed by an additive nonlinear  PDE. In particular the density $\rho$ is a solution to the following first order PDE  
\[
\frac{\partial \rho}{\partial t} = -H_1 H_2 \frac{\partial}{\partial x}\left(D(t,x)\rho\right)
\]
where the diffusion coefficient satisfies 
\[
\frac{\partial D}{\partial t}+\frac{1}{2}\frac{\partial^2 D}{\partial x^2} = 0.
\]
 Note that in this case we do not require fractional operators for the equation governing $D$ and the PDE satisfied by the density $\rho$ is of the first order.

\section{Conclusions}
In this short note we recalled the definition of  fractional Brownian motion and illustrated some of its properties that led to the vast success of fBm in applications. Our main observation is the link between  fBm and a class of coupled PDEs featuring time-dependent diffusion coefficients and fractional derivatives/integrals. A similar system of PDE can be written for iterated fBms.   This representation may offer a new point of view for the study of fBm and its applications. Moreover one  could look into extensions of this representation to more general self-similar processes or processes exhibiting long-range dependence.

\end{document}